\documentclass[10pt]{article}
\usepackage{graphicx}
\usepackage{amsmath}

\usepackage{amstext}
\usepackage{amsthm}
\newtheorem{theorem}{Theorem}
\newtheorem{proposition}{Proposition}
\newtheorem{lemma}{Lemma}

\newtheorem{corollary}{Corollary}

\begin{document}

\begin{titlepage}

\vskip0.5truecm

\vskip1.0truecm

\begin{center}

{\LARGE \bf Mather's regions of instability for annulus diffeomorphisms}

\end{center}

\vskip  0.4truecm

\centerline {{\large Salvador Addas-Zanata and  F\'abio Armando Tal}}

\vskip 0.2truecm

\centerline { {\sl Instituto de Matem\'atica e Estat\'\i stica }}

\centerline {{\sl Universidade de S\~ao Paulo}}

\centerline {{\sl Rua do Mat\~ao 1010, Cidade Universit\'aria,}} 

\centerline {{\sl 05508-090 S\~ao Paulo, SP, Brazil}}

\vskip 0.7truecm

\begin{abstract}

Let $f$ be a $C^{1+\varepsilon}$ diffeomorphism of the closed annulus $A$ that preserves 
orientation and the boundary components, and $\widetilde{f}$ be a lift of $f$ to its universal covering space. Assume that $A$ is a Birkhoff region of instability
 for $f$, and the rotation set of $\widetilde{f}$ is a non-degenerate interval. Then there exists an open $f$-invariant essential annulus $A^*$ whose frontier 
intersects both boundary components 
of $A$, and points $z^+$ and $z^-$ in $A^*$, such that the positive (resp. negative) orbit of $z^+$ converges to a set contained in the upper (resp. lower) 
boundary component of $A^*$ and the positive (resp. negative) orbit of $z^-$ converges to a set contained in the lower (resp. upper) 
boundary component of $A^*$. This extends a celebrated result originally proved by Mather in the context of area-preserving twist diffeomorphisms.

\end{abstract}

\vskip 0.3truecm

\vskip 2.0truecm



\vskip 0.8truecm

\noindent{\bf e-mail:} sazanata@ime.usp.br and fabiotal@ime.usp.br

\vskip 1.0truecm

%\noindent{\bf 2000 Mathematics Subject Classification:} 

\vfill
\hrule
\noindent{\footnotesize{Addas-Zanata would like to thank FAPESP (project
 2016/25053-8) and Tal would like to thank FAPESP (project
 2016/25053-8) and CNPq (project 303969/2018-0) for partial funding. }}

\end{titlepage}

\section{Introduction and statements of the results}

When studying the dynamics of $C^1$-area-preserving twist diffeomorphisms of
the closed annulus $A=S^1 \times [0,1],$ a celebrated theorem due to Mather \cite{mather}
states that if $f$ is such a diffeomorphism and $A$ contains no essential $f$%
-invariant continuum apart from each boundary component of $A$, except for 
continua containing both boundary components, (in other words, $A$ is
minimal with respect to the inclusion: its interior contains no proper 
essential $f$-invariant open annulus), then there are points $z^{+},z^{-}$ in $A$ such
that the $\alpha $-limit set of $z^{+}$ is contained in $S^1\times \{0\},$
the $\omega $-limit set of $z^{+}$ is contained in $S^1\times \{1\},$ the $%
\alpha $-limit set of $z^{-}$ is contained in $S^1\times \{1\}$ and the $%
\omega $-limit set of $z^{-}$ is contained in $S^1\times \{0\}.$ Mather's
theorem was proved using an intricate variational argument, which gave a lot
of insight about what happens in the $C^r$-generic situation, for all $r\geq
1.$ Later on, Le Calvez developed a completely topological proof \cite{inst}%
, extending Mather's result to twist area-preserving homeomorphisms.

In this paper, our main objective is to prove a version of this result in
the $C^{1+\varepsilon }$ world, without the twist and area-preservation
hypotheses and under the weaker condition of the annulus being a Birkhoff
region of instability, a situation which was widely considered for twist
homeomorphisms, see for instance \cite{birk}, \cite{mather}, \cite{inst} and 
\cite{frankslecal}. Namely, we consider $C^{1+\varepsilon }$ diffeomorphisms
$f:A\rightarrow A$ which preserve orientation and the boundary components,
whose rotation sets are non-degenerate intervals and as explained above, we
assume that $A$ is a Birkhoff region of instability for $f$ (see below for
the precise definition). The non-degeneracy condition on the rotation set is
always satisfied by twist maps. 

One last remark is that, in the area-preserving world, if $A$ is minimal as
explained above, then it is also a Birkhoff region of instability (this was
originally proved by Birkhoff for twist maps in \cite{birk}), but the two
definitions are not equivalent.

In order to state our main theorems properly, below we present some
definitions.

\vskip0.2truecm

{\bf Definitions.}

\begin{enumerate}
\item  Let $\rm{Diff}_0^r(A)$ be the subset of $C^r$ (for any $r\geq 0)$
diffeomorphisms $f:A\rightarrow A$ (when $r=0,$ $f$ is just a homeomorphism)
which preserve orientation and the boundary components of $A=S^1\times [0,1].
$ A lift of $f$ to the universal cover of the annulus $\widetilde{A}=${\rm I}%
\negthinspace 
{\rm R}$\times [0,1],$ is denoted $\widetilde{f}${\rm $,$} a homeomorphism
which satisfies $\widetilde{f}(\widetilde{z}+(1,0))=\widetilde{f}(\widetilde{%
z})+(1,0)$ for all $\widetilde{z}\in \widetilde{A}$.

\item  When $r=1+\varepsilon $ for some $0<\varepsilon <1$ in the above definition,
we mean that $Df$ is $\varepsilon $-Holder.

\item  The annulus $A$ is said to be a {\it {Birkhoff region of instability}}
for some $f\in \rm{Diff}_0^0(A)$ if, for all $\varepsilon >0$, there exist
integers $N,M>0$ such that $f^N(S^1\times ]0,\varepsilon [)$ intersects
    $S^1\times ]1-\varepsilon ,1[$ and $f^{-M}(S^1\times ]0,\varepsilon [)$ also
intersects $S^1\times ]1-\varepsilon ,1[$. We say that $A$ is a {\it {Mather
region of instability}} for $f$ if there are points $z^{+},z^{-}$ in $A,$
such that the $\alpha $-limit set of $z^+$ is contained in
the lower boundary component of $A$ 
and its $\omega $-limit set is contained in the
upper boundary component of $A.$ Similarly, the $\alpha $-limit set of $z^-$ is contained in
the upper boundary component of $A$ and its $\omega $-limit set is contained in the
lower boundary component of $A$.

\item  Let $p_1:\widetilde{A}$ $\rightarrow ${\rm I}\negthinspace {\rm R} be
the projection on the horizontal coordinate and as usual, let $p:\widetilde{A%
}$ $\rightarrow A$ be the covering mapping. Fixed $f\in \rm{Diff}_0^0(A)$ and a
lift $\widetilde{f},$ the displacement function $\phi :A${\rm $\rightarrow $I%
}\negthinspace {\rm R} is defined as

$$
\phi (z)=p_1\circ \widetilde{f}(\widetilde{z})-p_1(\widetilde{z}), 
$$

for any $\widetilde{z}\in p^{-1}(z).$

\item  Given any $f\in \rm{Diff}_0^0(A)$ and fixed some lift $\widetilde{f},$ a
point $z\in A$ is said to have {\it rotation number $\rho _0$} if the limit $%
\lim _{n\to \infty }\frac 1n\sum_{i=0}^{n-1}\phi (f^i(z))$ exists and is
equal to $\rho _0$. We say that $\rho _0$ is {\it realized by the compact
set $K\subset A$} if $K$ is $f$-invariant, and all points in $K$ have
rotation number $\rho _0$.

\item  Given any $f\in \rm{Diff}_0^0(A)$ and fixed some lift $\widetilde{f},$ the
rotation set of $\widetilde{f}$ is defined as

$$
\rho (\widetilde{f})=\{\omega \in {\rm I}\negthinspace {\rm R:}\text{ there
exists a Borel probability }f\text{-invariant} 
$$
$$
\text{ measure }\mu \text{ such that }\omega =\int_A\phi (z)d\mu \}. 
$$

Clearly, from the convexity of the subset of Borel probability $f$-invariant
measures, $\rho (\widetilde{f})$ is a closed interval, maybe a single point.
Moreover, its extremes are realized by ergodic measures, see \cite{misiu}.
In this generality, not much more can be said. There are well-known examples
with a non-degenerate interval as a rotation set, for which only the
extremes are the rotation numbers of some orbits.

\item  Given $f\in \rm{Diff}_0^0(A),$ we say that it satisfies the curve
intersection property, if for any homotopically non-trivial simple closed
curve $\gamma \subset A,$ we have $f(\gamma )\cap \gamma \neq \emptyset .$
It is not hard to see that if some $f\in \rm{Diff}_0^0(A)$ satisfies the curve
intersection property, then $f^n$ also satisfies it, for all integers $n \neq 0$.
Also, it is immediate that, in case $A$ is a Birkhoff region of instability
for some $f\in \rm{Diff}_0^0(A),$ then $f$ satisfies the curve intersection
property.

\item  Let $f\in \rm{Diff}_0^0(A)$ and $\gamma $ contained in the interior
  of $A$
  %considered as a surface with boundary,
  be a homotopically non-trivial simple
closed curve. Denote by $\gamma ^{-}$ the connected component of $\gamma ^c$
which contains $S^1\times \{0\}$ and analogously, let $\gamma ^{+}$ be the
connected component of $\gamma ^c$ which contains $S^1\times \{1\}.$ Denote
by $\xi^-(\gamma )$ the connected component of the
maximal invariant set in the closure of $\gamma ^-$ which contains
$S^1\times \{0\}$ and by $\xi^+(\gamma )$ the connected component of the
maximal invariant set in the closure of $\gamma ^+$ which contains
$S^1\times \{1\}.$
Finally, for all integers $n>1,$
denote by $\xi $$_{1/n}^{-}=$ $\xi ^{-}(S^1\times \{1-1/n\})$ and $\xi $$%
_{1/n}^{+}=$ $\xi $$^{+}(S^1\times \{1/n\}).$
\end{enumerate}

Note that, in case $A$ is a Birkhoff region of instability, 
$$
\partial (\xi ^{-}(\gamma ))^c\cap S^1\times \{0\}\neq \emptyset \text{ and }%
\partial (\xi ^{+}(\gamma ))^c\cap S^1\times \{1\}\neq \emptyset . 
$$

We are ready to state our main theorem.

\begin{theorem}
  \label{um} Let $f\in \rm{Diff}_0^{1+\varepsilon}(A)$ for some $\varepsilon>0$ be
  such that $A$ is a
Birkhoff region of instability and for some fixed lift $\widetilde{f},$ $%
\rho (\widetilde{f})$ has interior. Then there exists a homotopically
non-trivial simple closed curve $\gamma \subset A$ and an $f$-invariant
minimal open annulus $A^{*}\subset A$ containing $\gamma $ such that $%
\partial A^{*}$ intersects both $S^1\times \{0\}$ and $S^1\times \{1\}$ and $%
A^{*}$ is a Mather region of instability.
\end{theorem}

{\bf Remark:} Clearly, if $A$ is itself minimal as in the hypothesis of
Mather's original theorem, then $A=A^{*}$ and so, the whole annulus is a
Mather region of instability.

\vskip0.2truecm

Note that, in Proposition 2.19 of \cite{contal}, a result in the same
direction was obtained. There, $f$ was an area-preserving homeomorphim,
instead of a $C^{1+\varepsilon}$ diffeomorphim, and the thesis obtained was that $%
A^{*} $ is a mixed $SN$ region of instability, a condition weaker than being
a Mather region of instability.

Let us comment on the hypothesis and thesis of Theorem~\ref{um}. First, we
point out that the result does not hold in case the rotation set of $
\widetilde{f}$ is a single point. There are known examples (see for instance 
\cite{barney}) of smooth maps $f:A\rightarrow A$ having a single rotation
number that are both weak-mixing (therefore $A$ is a Birkhoff region of
instability), but also rigid (meaning that there is a sequence of positive
iterates of $f$ that converges to the identity). So, no subannulus of $%
f$ can be a Mather region of instability.

\begin{figure}[ht!]
\hfill
\includegraphics [height=48mm]{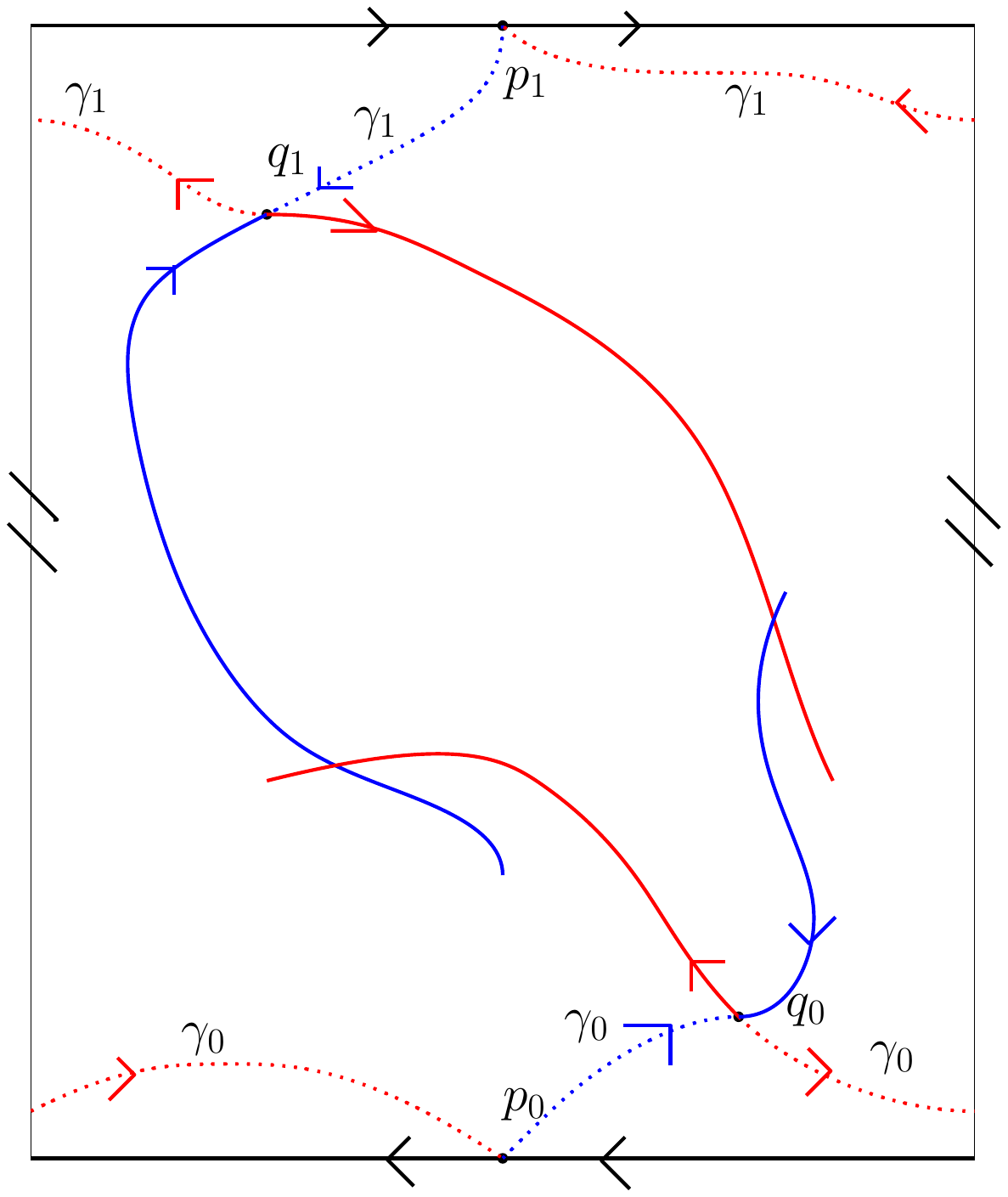}
\hfill{}
\caption{\small Sketch of the smooth diffeomorphism for which $A$ (obtained by gluing the two lateral sides together) is not a Mather Region of Instability.}
\label{figure_example}
\end{figure}

A natural question is also to understand, under our hypotheses, if $A$
itself is always a Mather region of instability. But this is false, and we sketch an example of a $C^{\infty}$ area-preserving diffeomorphism. Take $f:A\rightarrow A$, which extends to a smooth area-preserving diffeomorphism $g:S^1\times${\rm I}%
\negthinspace 
{\rm R}$\rightarrow S^1\times${\rm I}%
\negthinspace 
{\rm R} such that the restriction of $f$ to the upper boundary has a single degenerate topological
saddle fixed point $p_1$, that is, $p_1$ is a fixed point such that the differential of $g$ at $p_1$ is the identity, and such that there exists a local topological conjugation between the dynamics of $g$ at $p_1$ and a linear hyperbolic saddle at the origin. Likewise, we assume that $f$ has a single degenerate topological saddle fixed point $p_0$ in the lower boundary. So that each boundary of $A$ consists of a single fixed point and a saddle connection.  Furthermore, we assume that $f$ has two hyperbolic saddle points, $q_0$ and $q_1$, and that there exists saddle connections between a branch of the unstable manifold of $q_0$ and a branch of the stable manifold of $p_0$, as well as a connection between a branch of the stable manifold of $q_0$ and a branch of the unstable manifold of $p_0$, so that there exists an essential invariant closed curve $\gamma_0$ that intersects the lower boundary just at the point $p_0$. One can do a similar picture, with two saddle connections between $p_1$ and $q_1$, and an invariant essential closed curve $\gamma_1$ made by these connections and the points $p_1$ and $q_1$, which intersects the upper boundary just at $p_1$. Finally, one can assume that there are transversal heteroclinic intersections between the ``free'' branches of the unstable manifold of $q_0$ and the stable manifold of $q_1$, and between  the unstable manifold of $q_1$ and the stable manifold of $q_0$, see Figure~\ref{figure_example}. Finally, we may assume that, for a given lift $\widetilde f$ of $f$, $p_0$ and $p_1$ have different rotation numbers.

In this picture, from the $\lambda$-lemma one would get that the unstable manifold of $q_1$ accumulates on $p_0$ and since the unstable manifold of $q_1$ is contained in the closure of the future orbit of any neighborhood of $p_1$, one gets that there exists points arbitrarily close to $p_1$ whose future orbit lie arbitrarily close to $p_0$. A similar argument shows that there exists points arbitrarily close to $p_0$ whose future orbit lie arbitrarily close to $p_1$, and so  $A$ is a Birkhoff region of instability. But we claim $A$ cannot be a Mather region of instability for $f$. Indeed, if the $\omega$-limit of a point $z$ in the included in $S^1\times \{1\}$, then either $z$ lies in a stable branch of $p_1$, in which case the $\alpha$-limit  of $z$ is either $\{p_1\}$ or $\{q_1\}$, or $z$ must lie above the graph determined by $\gamma_1\cup\{p_1\}$. Since the closure of the later region is invariant and disjoint from $S^1\times \{0\}$, this implies that the $\alpha$-limit of $z$ is disjoint from $S^1\times \{0\}$. This shows that there does not exist a point whose $\alpha$-limit lies in $S^1\times \{0\}$ and whose $\omega$-limit lies in  $S^1\times \{1\}$.

%We believe, but have no current example, that using the $C^1$ perturbation method developed by Bonatti and Crovisier in \cite{bonacro}, it may be possible to assume that the above $f$
%is a $C^1$-area-preserving twist diffeomorphism. If this holds, as $f$ is transitive in the
%essential connected component of $\gamma ^c$, such example would show that even for area-preserving twist
%diffeomorphisms, $A$ being a Birkhoff region of instability does not imply
%that it is a Mather region of instability.

The next lemma is a crucial step in the proof of Theorem~\ref{um}. It also
explains how the annulus $A^{*}$ is constructed.

\begin{lemma}
\label{umaux} Under the hypotheses of Theorem 1, there exists a
homotopically non-trivial closed curve $\gamma$ contained in the interior of $A$
such that $\xi^-(\beta )\cap \gamma =\emptyset $ and 
$\xi^+(\beta )\cap \gamma =\emptyset $  for any
homotopically non-trivial simple closed curve $\beta \subset A.$ So, there
exist $f$-invariant continua $K^{-}\supset S^1\times \{0\}$ and $%
K^{+}\supset S^1\times \{1\}$ such that $K^{-}\cap K^{+}=\emptyset $ and for
all sufficiently large integers $n>1,$ $\xi $$_{1/n}^{-}=K^{-}$ and $\xi $$%
_{1/n}^{+}=K^{+}.$ Moreover, if $A^{*}$ is the $f$-invariant open annulus
between $K^{-}$ and $K^{+},$ then $\partial A^{*}$ intersects both $%
S^1\times \{0\}$ and $S^1\times \{1\}.$
\end{lemma}

%For area-preserving twist maps 
When $A$ is minimal, in the sense that its interior contains no $f$-invariant proper
essential open sub-annulus, then $\xi $$_{1/n}^{-}=S^1\times \{0\}$ and $\xi 
$$_{1/n}^{+}=S^1\times \{1\},$ for all integers $n>1,$ but when $A$ is just
a Birkhoff region of instability, we cannot avoid considering the sets $%
K^{-} $ and $K^{+}$ from Lemma \ref{umaux}.

We also provide the following result, which is fundamental in the proof of
Theorem \ref{um}, but whose interest stands alone for its possible
applications:

\begin{theorem}
\label{dois} Again, under the hypotheses of Theorem 1, there exists $E,$ an
open and dense subset of $\rho (\widetilde{f}),$ such that for any rational
number $p/q\in E,$ there exists a hyperbolic periodic saddle point $z$
contained in the interior of $A$ whose rotation number is $p/q$ and a
homotopically non-trivial closed curve $\gamma _{p/q}\ni z,$ contained in
the union of the stable and unstable manifolds of $z.$
\end{theorem}

Finally, we remark that part of the interest in regions of instability comes
from the fact that the dynamics in these regions is usually very rich. For
instance, a classical result says that, if $f\in \rm{Diff}_0^0(A)$ has the curve
intersection property (a condition that is satisfied when $A$ is a Birkhoff
region of instability), and if $\rho (\widetilde{f})=[a,b]$, then for any
rational number $a<p/q<b,$ $\widetilde{f}^q-(p,0)$ has a fixed point. In
other words, all rationals in the interior of the rotation set are realized
by periodic orbits. In order to prove this statement, note that in \cite
{misiu} it was proved that $a$ and $b$ are realized by ergodic Borel $f$%
-invariant probability measures, that is, $a$ and $b$ are equal to the
rotation numbers of actual points in $A$. So, if for some $a<p/q<b$, $
\widetilde{f}^q-(p,0)$ does not have a fixed point, then the version of
Brouwer's translation theorem applied to the annulus that appears in a
1928-29 paper of Kerekjarto \cite{kerek} implies the existence of a
homotopically non-trivial simple closed curve disjoint from its image under $%
f$, a contradiction with the curve intersection property.

Moreover, if the twist condition is present, even for any irrational number $%
\rho_0$ in the rotation set, one can find an $f$-invariant compact set $%
K_{\rho_0}$ that realizes $\rho_0$, such that the restriction of $f$ to $%
K_{\rho_0}$ is semi-conjugated to the irrational rotation of the circle with
the same rotation number. The fact that the full rotation set is realized by
compact $f$-invariant sets was also proved in the absence of the twist
condition, for area preserving homeomorphims (\cite{lecam, contal}).

Our final result, a direct consequence of Theorems~\ref{um}, \ref{dois} and
Theorem C of \cite{contal}, says that:

\begin{theorem}
  \label{tres} Let $f\in \rm{Diff}_0^{1+\varepsilon}(A)$ for some $\varepsilon>0$
  be such that $A$ is a
Birkhoff region of instability and let $\widetilde{f}$ be a lift of $f$ to
its universal covering space. Then there are at most two numbers in $\rho (
\widetilde{f})$ that are not realized by compact $f$-invariant sets.
\end{theorem}

The paper is organized as follows. In the next section we present all the
necessary preliminary results, with a brief overview and remainder of the
tools needed in this work. Section~3 is dedicated to the proofs of our main
results.

%, the same holds for all its iterates. And the following theorem
%holds:

%\begin{theorem}
%\label{rotset} Let $f\in Diff_0^0(A)$ satisfies the C.I.P. and suppose
%that $\rho (\widetilde{f})=[a,b]$ for some $a<b.$ Then for all $a<p/q<b,$ $
%\widetilde{f}^q-(p,0)$ has a fixed point. In other words, all rationals in
%the interior of the rotation set are realized by periodic orbits.
%\end{theorem}

%Probably the first place where a proof of this result appeared is in a
%1928-29 paper of Kerekjarto \cite{kerek}. (bla maybe Birk second paper???)

%Clearly, when $A$ is a Birkhoff region of instability, the C.I.P. holds.
%So the above theorem implies that for $f\in Diff_0^0(A)$ for which $A$ is a
%Birkhoff region of instability, all rational interior points of the
%rotation set are realized by periodic orbits.
%Our main theorems are proved in reverse order.

\section{Preliminaries}

In this section we describe some theories we use and quote some results.

\subsection{Prime ends compactification of open disks}

%\label{sec:primends}

If $D$ is an open topological disk of an oriented surface $S$ potentially
with boundary, such that $\partial D$ is a Jordan curve and $f$ is an
orientation preserving homeomorphism of that surface which satisfies $%
f(D)=D, $ then $f:\partial D\rightarrow \partial D$ is conjugate to a
homeomorphism of the circle, and so a real number $\rho (D),$ the rotation
number of $f\mid _{\partial D}$ can be associated to this problem. By the
classical properties of rotation numbers, if $\rho (D)$ is rational, then
there exists a periodic point in $\partial D$ and if it is not, then there
are no such points. This is known since Poincar\'e. The difficulties arise
when we do not assume $\partial D$ to be a Jordan curve.

The prime ends compactification is a way to attach to $D$ a circle called
the circle of prime ends of $D,$ obtaining a space $D\sqcup S^1$ with a
topology that makes it homeomorphic to the closed unit disk. If, as above,
we assume the existence of an orientation preserving homeomorphism $f$ of $S$
such that $f(D)=D,$ then $f\mid _D$ extends to $D\sqcup S^1.$ The prime ends
rotation number of $f$ in $D,$ still denoted $\rho (D),$ is the usual
rotation number of the orientation preserving homeomorphism induced on $S^1$
by the extension of $f\mid _D.$ But things may be quite different in this
setting. In full generality, it is not true that when $\rho (D)$ is
rational, there are periodic points in $\partial D$ and for some examples, $%
\rho (D)$ is irrational and $\partial D$ is not periodic point free. Anyway,
the only result on this subject we need is the following classical lemma (as
usual, a point $z\in \partial D$ is said to be accessible if there exists a
simple arc $\gamma :[0,1]\rightarrow D\cup \partial D$ such that $\gamma
([0,1[)\subset D$ and $\gamma (1)=z)$ whose proof, for instance, can be found in
Theorem 16 of \cite{ortega}.

\begin{lemma}
\label{acessperiod} Let $f\in \rm{Diff}_0^0(A)$ and let $D\subset A$ be an $f$%
-invariant open annulus given by the complement of some $f$-invariant
continuum $K$ which contains $S^1\times \{0\}$ and avoids $S^1\times \{1\}.$
Then the boundary of $D$ has two connected components, one is $S^1\times
\{1\}$ and the other one is some continuum $M\subseteq K.$ Moreover, if $%
z_1,z_2\in M$ are periodic points, both accessible from $D,$ then they have
the same rotation number.
\end{lemma}

{\bf Remark: }As $D$ is not a disk, by prime ends rotation number of $D$
(denoted $\rho (D)),$ we mean the following: Contract $S^1\times \{1\}$ to a
point $N$ %and $S^1\times \{0\}$ to a point $S$ 
in order to turn $A$ into a closed disk. Clearly, $f:A\rightarrow A$ induces
a homeomorphism of this closed disk which fixes $N.$ Now, $D$ becomes an
open topological disk and $\rho (D)$ is the prime ends rotation number of
this disk.

\vskip0.2truecm

Some interesting facts related to the previous result, but not necessary
for us are the following: $\rho (D)$ is equal to the rotation number (in $A$)
of any accessible periodic point in $M.$ More generally, the main result of \cite{ortega}
says that for any accessible point in $M=\partial D,$ either its
forward or backward annulus rotation number is equal to $\rho (D).$ 

For more information on the theory of prime ends, see for instance \cite{mather}, \cite
{koropatmey} and \cite{ortega}.

\subsection{Topologically transverse intersections}

Let $S$ be a compact orientable surface. We say that a closed topological
disk $R\subset S$ is a topological rectangle if its boundary, which is a
Jordan curve, is given by the union of four $C^{1}$ oriented arcs, $\gamma
_{1},\gamma _{2},\gamma _{3},\gamma _{4}$ such that: the end point of $%
\gamma _{1}$ is the first point of $\gamma _{2}$, the end point of $\gamma
_{2}$ is the first point of $\gamma _{3}$, the end point of $\gamma _{3}$ is
the first point of $\gamma _{4}$ and the end point of $\gamma _{4}$ is the
first point of $\gamma _{1}$. We also assume that for $i,j\in \{1,2,3,4\},$ $%
i\neq j,$ the intersection between $\gamma _{i}$ and $\gamma _{j}$ is either
empty or $C^{1}$-transversal.

Suppose $f$ is a diffeomorphism of $S$ and assume it has a hyperbolic $n$%
-periodic saddle point $p$ contained in the interior of $S.$

\vskip0.2truecm

{\bf Definition of topological transversality:} We say that some continuum $%
K\subset S$ has a topologically transverse intersection with a branch $%
\alpha $ at $p$ (stable or unstable), if there exists a topological
rectangle $R\subset S$, whose boundary is given by arcs $\gamma _{1},\gamma
_{2},\gamma _{3},\gamma _{4}$ as explained above, and there exists an arc $%
\alpha ^{\prime }\subset \alpha $, whose interior is contained in the
interior of $R$ and its extreme points belong, one to the interior of $%
\gamma _{1}$ and the other one to the interior of $\gamma _{3}$ (the
intersections between $\gamma _{1}$ and $\alpha ^{\prime }$ and $\gamma _{3}$
and $\alpha ^{\prime } $ are both $C^{1}$-transverse), such that $K$
contains a subcontinuum $K^{\ast }\subset R$ which intersects both $\gamma
_{2}$ and $\gamma _{4}$, and avoids $\gamma _{1}$ and $\gamma _{3}$.

\vskip0.2truecm

Such intersections are important because of the following result:

\begin{proposition}
\label{toptrans} In the above setting, if $K$ is a continuum which has a
topologically transverse intersection with a stable branch $\alpha $ at $p$
(the unstable case is analogous), then given any $\mu >0$ and a compact arc $%
\theta \subset W^u(p),$ there exists $M>0$ such that for all $m\geq M,$ $%
f^{m.2n}(K)$ contains a continuum $\theta _m$ which is $\mu $-close to $%
\theta $ for the Hausdorff distance.
\end{proposition}

\textit{Proof. }By considering $f^{2n}$ instead of $f,$ we can assume that $%
p $ is fixed and all branches at $p$ are $f$-invariant. From the
Hartman-Grobman theorem, there exists $V,$ an open neighborhood of $p,$ $%
W\subset \mathrm{I}\!\mathrm{R}^{2},$ an open neighborhood of the origin and
a homeomorphism $\varphi :V\rightarrow W$ which conjugates $f\mid _{V}$ to
the linear model $H(x,y)=(x/2,2y)$ restricted to $W$. Suppose, without loss
of generality, that, locally, $\varphi(\alpha)$ corresponds to the positive $%
x$-axis in $W.$

Consider the rectangle $R_{c ,d }=[-c,c]\times \lbrack -d,d]$ for $c ,d >0$
such that $H(R_{c ,d })\cup R_{c,d}\subset W.$ As $K$ is topologically
transverse to $\alpha$, there exists a topological rectangle $R,$ such that $%
K^{\ast }\subset R$ is a subcontinuum of $K$ which intersects both $\gamma
_{2}$ and $\gamma _{4},$ and avoids $\gamma _{1}$ and $\gamma _{3}$ (see the
definition of topological transversality above). It is clearly possible to
modify $R$ by choosing $\gamma _{2}^*$ and $\gamma _{4}^*$ much closer (in a 
$C^{1}$-way) to $\alpha ^{\prime }$ (the connected arc contained in $\alpha$ whose interior
is contained in $R$ and whose extreme points are contained, one in the interior
of $\gamma_1$ and the other, in the interior of $\gamma_3$), 
and choosing $\gamma_1^*$ and $%
\gamma_3^*$ subarcs of $\gamma_1$ and $\gamma_3$ respectively, in a way that 
$\gamma_1^*, \gamma_2^*, \gamma_3^*, \gamma_4^*$ form the boundary of a new
rectangle $R^*$, such that for some integer $N\geq 0,$ $f^{N}(R^*)\subset V$
and the corresponding rectangle $R^{\prime}=\varphi (f^{N}(R^*))$ belongs to 
$]0,c \lbrack \times ]-d ,d\lbrack \subset R_{c,d}$. Furthermore, there
exists some subcontinuum $K^{**}$ of $K^*$ which is disjoint from $%
\gamma_1^* $ and $\gamma_3^*$, intersects both $\gamma_2^*$ and $\gamma_4^*$
and is contained in $R^*$. Related to $K^{**}$, let us choose real numbers $%
a,b,\delta >0$ such that $0<a<b<c ,$ $0<\delta <d ,$ $R^{\prime}\subset
]a,b[\times ]-d,d[$ and $[-c,c]\times \lbrack-\delta,\delta]\cap \varphi
(f^{N}(\gamma _{2}^*\cup \gamma _{4}^*))=\emptyset.$ Clearly, both $%
[a,b]\times\{\delta\}$ and $[a,b]\times\{-\delta\}$ intersect $R^{\prime}$.

Let $\Gamma \subset [a,b]\times [-\delta,\delta]$ be a connected component
of 
$$
[a,b]\times [-\delta,\delta]\cap \varphi (f^N(K^{**})) 
$$
which intersects both $[a,b]\times \{-\delta\}$ and $[a,b]\times \{\delta\}.$
Clearly, $\Gamma\cap\{a,b\}\times [-\delta,\delta ]=\emptyset .$

\vskip0.2truecm %Now, we claim:
{\bf \noindent Claim.} {\it For any $\varepsilon>0$, there exists $M>0$ such
that for all integers $m\geq M,$ there is a subcontinuum $\Gamma^{\prime
}\subset \Gamma$, depending on $\varepsilon$ and $m$, such that $%
H^{m}(\Gamma^{\prime })\subset \lbrack0,\varepsilon \lbrack \times \lbrack
-d ,d ]$ and $H^{m}(\Gamma ^{\prime })$ intersects both $[0,\varepsilon
\lbrack \times \{-d \}$ and $[0,\varepsilon \lbrack \times \{ d \}.$ }

\vskip0.2truecm

Let us, before proving the claim, show that it implies the proposition.
Indeed, if it holds, then as $\varepsilon \rightarrow 0,$ $\varphi
^{-1}(H^{m}(\Gamma^{\prime }))\subset f^{m+N}(K^{\ast \ast})$ converges in
the Hausdorff topology to $\varphi ^{-1}(\{0\}\times \lbrack -d ,d ]),$ a
local unstable manifold at $p. $

For $\theta$ as in the statement of the proposition, there exists some
integer $J>0$ such that $\theta\subset f^J(\varphi ^{-1}(\{0\}\times [-d,d
]))$. As $\varepsilon \rightarrow 0,$ $f^J\left( \varphi
^{-1}(H^m(\Gamma^{\prime }))\right) \subset f^{J+m+N}(K^{**})$ converges in
the Hausdorff topology to $f^J\left( \varphi ^{-1}(\{0\}\times [-d ,d
])\right),$ which contains $\theta$. So, there exists a subcontinuum $\theta
_m$ of $f^{J+m+N}(K^{**})$ that converges to $\theta $ in the Hausdorff
topology as $\varepsilon \rightarrow 0.$

\vskip0.2truecm

Therefore, to conclude the proof of Proposition \ref{toptrans}, we have to
show that the above claim holds. For this, given $\varepsilon>0$, let $M>0$
be an integer such that:

\begin{itemize}
\item  $b/2^M<\varepsilon;$

\item  $d/2^M<\delta ,$
\end{itemize}

Clearly, $M\rightarrow \infty $ as $\varepsilon \rightarrow 0.$ Moreover,
for all $m\geq M,$ let $\Gamma^{\prime }$ be a connected component of $%
\Gamma \cap \lbrack a,b]\times \lbrack -d /2^{m},d /2^{m}]$ which intersects
both $[a,b]\times \{-d /2^{m}\}$ and $[a,b]\times \{d /2^{m}\}.$ From the
choice of $M>0,$ $\Gamma^{\prime }$ is not empty and it clearly satisfies $%
H^{m}(\Gamma^{\prime })\subset \lbrack 0,\varepsilon \lbrack \times \lbrack
-d ,d ]$ and $H^{m}(\Gamma^{\prime })$ intersects both $[0,\varepsilon
\lbrack \times \{-d \}$ and $[0,\varepsilon \lbrack \times \{d \}.$

This proves the claim and concludes the proof of the proposition. $\qed$

\vskip0.2truecm

Still in the above setting, if $\theta,$ the compact subarc of a branch, 
has a topologically transverse intersection with some other continuum $T,$
then $f^{m.2n}(K)$ also intersects $T$ provided $m>0$ is large enough. For
more about topologically transverse intersections, see for instance \cite
{c1epsilon} and \cite{bruno}.

The next result will be used in the proof of Theorem 1.

\begin{lemma}
\label{l.pertseparada} In the above setting, assume $K$ is a continuum which
does not have a topologically transverse intersection with the stable branch 
$\alpha $ at a saddle $p$. Then, for every $\varepsilon>0$ and any $\theta,$
subarc of $\alpha$ such that $K$ is disjoint from the endpoints of $\theta ,$
there
exists a simple arc $\theta _{\varepsilon}$ that is contained in the
$\varepsilon$-neighborhood of $\theta $, has the same endpoints as $\theta $, and is
disjoint from $K$.
\end{lemma}

{\it Proof.} Fix $\varepsilon>0$ and $\theta$ as in the statement of the
lemma. Let $f_\theta(t),$ $t\in[0,1]$ be a parametrization of $\theta.$ One
can find an $\varepsilon$-neighborhood $V$ of $\theta$, a neighborhood $W$
of $[0,1]\times\{0\}$ and a $C^1$ diffeomorphism $\phi:V\to W$ such that $%
\phi(f_\theta(t))= (t,0)$. Also, if $\delta>0$ is sufficiently small, then $%
R= [0,1]\times[-\delta,\delta]$ is a subset of $W$ and $\phi(K \cap V)$ is
disjoint from $\{0,1\}\times[-\delta,\delta]$ and not contained in $R$. 
%$V_R$, defined as the union of the two vertical sides of $R$.%, $\{0\}\times[-\delta,\delta]$ and $\{1\}\times[-\delta,\delta].$
Consider the subset of $R,$ denoted $F:= \phi(K \cap V)\cap R.$ If it
separates $\{0\}\times[-\delta,\delta]$ from $\{1\}\times[-\delta,\delta],$
then, as $F$ is closed, there is a connected component of $F$ that separates
the former two sets in $R$. So, this component intersects both $%
[0,1]\times\{-\delta\}$ and $[0,1]\times\{\delta\},$ something that
contradicts the assumption that $K$ does not have a topologically transverse
intersection with $\alpha$.

In this way, there is a connected component $B$ of $R\setminus F$ which
contains $\{0\} 
\times[-\delta,\delta]$ and $\{1\}\times[-\delta,\delta]$. And it is open. So if we pick $\beta:[0,1]\to B$, a simple arc joining $(0,0)$ and $(1,0),$ it suffices to take $\theta_{%
\varepsilon}=\phi^{-1}(\beta)$. \qed
%
%
%there is acontains  $\{0\}\times[-\delta,\delta]$  also contains $\{1\}\times[-\delta,\delta]$. And it is open. So if we pick $\beta:[0,1]\to R$, a simple arc joining $(0,0)$ and $(1,0),$ contained in $B$, it suffices to take $\theta_{\varepsilon}=\phi^{-1}(\beta)$. \qed
%
%\vskip0.2truecm
%$F$ that separates the former two sets.  
%The fact that $K$ does not have a topologically transverse intersection with $\theta$ implies that every connected component of $F= \phi(K \cap V)\cap R$ is either contained in $(0,1)\times[-\delta, 0]$ or in $(0,1)\times[0,\delta]$, so each connected component of $F$ intersects either $[0,1]\times\{-\delta\}$ or $[0,1]\times\{\delta\}$. Denote by $F^-$ the union of the conneted components of $F$ that intersect the former, and by $F^+$ to the union of the connected components of $F$ that intersect the later. Let also $0<\kappa<\delta$ be smaller than
%$$
%\frac{1}{4}min\{distance(F^+,F^-),distance(F, \{0,1\}\times[-\delta,\delta])\}
%$$
%smaller than the distance between $F_1$ and $F_2$. 
%One verifies that the connected component $B$ of $R\setminus\left(F \cup \overline{B_{\kappa}(F^-)\cup B_{\kappa}(F^+)}\right)$ that contains  $\{0\}\times[-\delta,\delta]$  also contains $\{1\}\times[-\delta,\delta]$. And it is open. So if we pick $\beta:[0,1]\to R$, a simple arc joining $(0,0)$ and $(1,0),$ contained in $B$, it suffices to take $\theta_{\varepsilon}=\phi^{-1}(\beta)$. \qed
%
\vskip0.2truecm

\subsection{Some Pesin theory}

\label{pesin}

In this subsection, assume that $f:A\to A$ is a $C^{1+\varepsilon}$ diffeomorphism,
for some $\varepsilon>0.$ Recall that an $f$-invariant Borel probability measure $\mu $
is hyperbolic if all the Lyapunov exponents of $f$ are non-zero at $\mu $%
-almost every point (for instance, see the supplement of \cite{akbh}).
Remember that for $\mu $-almost every $z\in A,$ there are two Lyapunov
exponents $\lambda _{+}(z)\geq \lambda _{-}(z)$ defined as follows:

$$
\lambda _{+}(z)=\lim _{n\rightarrow \infty }\frac 1n \log \left\| Df^n(z)\right\| 
\text{ and }\lambda _{-}(z)=-\lim _{n\rightarrow \infty }\frac 1n \log \left\|
Df^{-n}(z)\right\| \text{ } 
$$

The next paragraphs were taken from \cite{asmp03}. They consist of an
informal description of the theory of non-uniformly hyperbolic systems,
together with some definitions and lemmas from \cite{asmp03}.

Let $\mu $ be a non-atomic hyperbolic ergodic $f$-invariant Borel
probability measure. Given $0<\delta <1,$ there exists a compact set $%
\Lambda _\delta $ (called Pesin set) with $\mu (\Lambda _\delta )>1-\delta ,$
having the following properties: for every $p\in \Lambda _\delta ,$ there
exists an open neighborhood $U_p,$ a compact neighborhood $V_p\subset U_p$
and a diffeomorphism $F:(-1,1)^2\to U_p,$ with $F(0,0)=p$ and $%
F([-1/10,1/10]^2)=V_p,$ such that:
\begin{itemize}

\item The local unstable manifolds $W_{\mathrm{loc}}^u(q)$, 
$q \in \Lambda _\delta \cap V_p$, given by the connected component of the set of $z\in
U_p$ sucht that dist$(f^n(z),f^n(q))\rightarrow 0$ as $n\rightarrow -\infty$ that
contains $q$, are the
images under $F$ of graphs of the form $\{(x,F_2(x)):x\in (-1,1)\},$ $F_2$ a
function with $k$-Lipschitz constant, for some $0<k<1.$ Any two such local
unstable manifolds are either disjoint or equal and they depend continuously
(in the Hausdorff topology) on the point $q\in \Lambda _\delta \cap V_p.$
\item Similarly, local stable manifolds $W_{\mathrm{loc}}^s(q)$, 
$q \in \Lambda _\delta \cap V_p$, given by the connected component of the set of $z\in
U_p$ such that dist$(f^n(z),f^n(q))\rightarrow 0$ as $n\rightarrow \infty$ that
contains $q$, are the images under $F$ of graphs of the form $%
\{(F_1(y),y):y\in (-1,1)\}$, $F_1$ a function with $k$-Lipschitz constant,
for some $0<k<1$. Any two such local stable manifolds are either disjoint or
equal and they depend continuously (again in the Hausdorff topology) on the
point $q\in \Lambda _\delta \cap V_p.$ 
\end{itemize}
These are the properties that
characterize a Pesin set.

It follows that there exists a continuous product structure in $\Lambda
_\delta \cap V_p:$ given any $r,r^{\prime }\in \Lambda _\delta \cap V_p,$ the
intersection $W_{loc}^u(r)\cap W_{loc}^s(r^{\prime })$ is transversal and
consists of exactly one point, which will be denoted $[r,r^{\prime }]$. This
intersection varies continuously with the two points and may not be in $%
\Lambda _\delta .$ Hence we can define maps $P_p^s:\Lambda _\delta \cap
V_p\to W_{loc}^s(p)$ and $P_p^u:\Lambda _\delta \cap V_p\to W_{loc}^u(p)$ as $%
P_p^s(q)=[q,p]$ and $P_p^u(q)=[p,q].$

Let $R^{\pm }$ denote the set of all points in $A$ which are both forward and
backward recurrent. By the Poincar\'e recurrence Theorem, $\mu (R^{\pm })$
is equal to $1.$

\begin{description}
\item[Definition]  {\bf (Accessible and inaccessible points).} A point $p\in
\Lambda _\delta \cap V_p\cap R^{\pm }$ is inaccessible if it is accumulated
on both sides of $W_{loc}^s(p)$ by points in $P_p^s(\Lambda _\delta \cap
V_p\cap R^{\pm })$ and also accumulated on both sides of $W_{loc}^u(p)$ by
points in $P_p^u(\Lambda _\delta \cap V_p\cap R^{\pm })$. Otherwise, $p$ is
accessible.
\end{description}

After this definition, we can state two lemmas from \cite{asmp03} about
accessible and inaccessible points and the relation between these points and
nearby hyperbolic periodic points.

\begin{lemma}
\label{andre1} Let $q\in \Lambda _\delta \cap V_p\cap R^{\pm }$ be an
inaccessible point. Then there exist rectangles enclosing $q,$ having sides
along the invariant manifolds of two hyperbolic periodic saddle points in $V_p$ 
(which are corners of the rectangles) and having arbitrarily small diameter.
\end{lemma}

The boundary of such a rectangle is a Jordan curve made up of alternating
segments of stable and unstable manifolds, two of each. The segments forming
the boundary are its sides and the intersection points of the sides are the
corners. As explained above, two of the corners are hyperbolic periodic saddle points and 
the other two corners are $C^1$-transverse heteroclinic intersections.
A rectangle is said to enclose $p$ if its interior, which is an
open topological disk, contains $p.$

\begin{lemma}
\label{andre2} The subset of accessible points in $\Lambda _\delta \cap
V_p\cap R^{\pm }$ has $\mu $ measure equal to zero.
\end{lemma}

Another concept that will be a crucial hypothesis for us is positive
topological entropy. In the following we describe why.

When the topological entropy $h(f\vert_K)$ is positive, for some compact $%
f $-invariant set $K,$ by the variational principle, there exists an $f$%
-invariant Borel probability measure $\mu _0$ with $\rm{supp}(\mu _0)$ 
(the topological support of $\mu_0$) contained in $K$
and positive metric entropy $h_{\mu _0}(f).$ Using the ergodic decomposition of $\mu _0$ we find
an extremal point $\mu $ of the set of Borel probability $f$-invariant
measures, such that $\rm{supp}(\mu )$ is also contained in $K$ and $h_\mu (f)>0.$
Since the extremal points of this set are ergodic measures, $\mu $ is
ergodic. The ergodicity and the positiveness of the entropy imply that $\mu $
has no atoms and applying the Ruelle inequality (which, in our case, says
that $\lambda _{+}(z)=\lambda _{+}(\mu )\geq h_\mu (f)>0$ for $\mu $-almost
every $z\in A),$ we get that $\mu $ has a positive Lyapunov exponent, see 
\cite{akbh}. Working with $f^{-1}$ and using the fact that $h_\mu
(f^{-1})=h_\mu (f)>0,$ we see that $f^{-1}$ must also have a positive
Lyapunov exponent with respect to $\mu ,$ which is the opposite of the
negative Lyapunov exponent for $f$.

Hence when $K$ is a compact $f$-invariant set and the topological entropy of 
{\it f}$\mid_K$ is positive, there always exists an ergodic, non-atomic,
invariant measure supported in $K,$ with non-zero Lyapunov exponents, one
positive and one negative, the measure having positive entropy: That is, an
hyperbolic measure.

The existence of this kind of measure will be important for us because of
Lemmas \ref{andre1} and \ref{andre2}.

\subsection{Some forcing results}

This subsection is mostly based on the work of Le Calvez and the second
author \cite{lectal1, lectal2} on forcing theory for surface homeomorphims,
which we restrain to explain in more detail as not to substantially increase
the length of this paper. We refer the interested reader to the above works,
as well as \cite{contal}.

Given $f\in {\rm {Diff}_0^0(A)}$ and a lift $\widetilde f$, we say that $f$
has a {\it rotational topological horseshoe} if there exists, for some power 
$g=f^r $ of $f$:

\begin{itemize}
\item  A lift $\widetilde{g}=\widetilde{f}^r-(s,0)$ to $\widetilde{A}$,
where $s$ is an integer,

\item  A compact $g$-invariant subset $\Lambda \subset A$,

\item  A compact subset $\widetilde{\Lambda }\subset \widetilde{A}$ such
that $p(\widetilde{\Lambda })=\Lambda $ and the restriction of $p$ to $
\widetilde{\Lambda }$ is a homeomophism onto its image,

\item  An integer $M_0$, a compact metric space $Y$, a homeomorphim $T$ of $%
Y$, and a surjective continuous map $\pi _1:Y\to \Lambda $ semiconjugating $%
T $ and $g$, such that for each $x\in \Lambda $, the cardinality of $\pi
_1^{-1}(x)$ is not larger than $M_0$,

\item  A continuous surjective map $\pi _2:Y\to \Sigma _2:=\{0,1\}^{{\rm Z}%
\negthinspace \negthinspace {\rm Z}},$ semiconjugating $T$ to the shift map $%
\sigma :\Sigma _2\to \Sigma _2$ defined so that $(\sigma
(u))_j=(u)_{j+1},\,u\in \Sigma _2$,

\item  Also, if $y\in Y$, $x=\pi _1(y)$, $\tilde x=p^{-1}(x)\cap \widetilde{%
\Lambda }$, then $\widetilde{g}(\widetilde{x})\in \widetilde{\Lambda }$ if $%
(\pi _2(y))_0=0$ and $\widetilde{g}(\widetilde{x})\in \widetilde{\Lambda }%
+(1,0)$ if $(\pi _2(y))_0=1$.
\end{itemize}

In general terms, if $f$ has a rotational topological horseshoe, then modulo
some finite extension and taking a power of the dynamics, we obtain a
compact invariant set where the displacement of points in the lift can be
estimated by a symbolic coding, and for which every possible coding with $2$
symbols is admitted. We remark that it follows directly from the definition,
that the rotation set of $\widetilde g$ contains the interval $[0, 1]$, and
so the rotation set of $\widetilde f$ contains the interval $[s/r, (s+1)/r]$%
. We also say that the $\widetilde f$-rotation set of the horseshoe contains
the interval $[s/r, (s+1)/r]$. But we can obtain a little more, which will
be useful:

\begin{lemma}
\label{transitiveentropy} Let $f\in {\rm {Diff}_0^0(A)}$ and let $\widetilde{%
f}$ be a lift of $f$.  Let $\Lambda ,\widetilde{%
\Lambda },g,\widetilde{g},Y,T,\pi _1$ and $\pi _2$ be as above. Then there
exists a compact $g$-invariant set $\Lambda _0\subset \Lambda $, such that
the restriction of $g$ to $\Lambda _0$ is transitive and such that, for
every rational $0<p/q<1,$ one can find a $g^q$-invariant compact set $%
K_{p/q}\subset \Lambda _0$ satisfying:

\begin{enumerate}
\item  The restriction of $g^q$ to $K_{p/q}$ has strictly positive
topological entropy.

\item  If $\widetilde{K}_{p/q}=\pi ^{-1}(K_{p/q})\cap \widetilde{\Lambda }$,
then $\widetilde{g}^q(\widetilde{K}_{p/q})=\widetilde{K}_{p/q}+(p,0)$
\end{enumerate}
\end{lemma}

{\it Proof. }Let $u_0\in \Sigma_2$ be a point whose future $\sigma$ orbit is
dense, and for notation sake denote $L_0=\pi_2^{-1}(u_0)$. We claim that
there exists some $y_0\in L_0$ that is recurrent for $T$. Indeed, let us consider the set $\mathcal{T}$ of all closed and $T$-invariant subsets of $Y$ that have nonempty intersection with $L_0$, which is naturally ordered by inclusion, that is, where for $F_1, F_2\in \mathcal{T}$ we denote $F_1\preceq F_2$ if $F_1\subset F_2$. Note that $\mathcal{T}$ is not empty as $Y$ belongs to it.  Consider a chain $(F_i)_{i\in I}$ with each $F_i$ in $\mathcal{T}$ and such that for all $i,j$ in $I$, either $F_{i}\preceq F_{j}$ or $F_{j}\preceq F_{i}$. We claim that $ F=\bigcap_{i\in I}F_i$ also belongs to $\mathcal T$.  Indeed, $F$ is compact as it is the intersection of compact sets. We will show that $F\cap L_0$ is nonempty. If not, by compactness of $L_0$, since the complements of $F_i$ would form an open covering of $L_0$, there would be finitely many indices $\{i_1, i_2, ..., i_n\}$ such that $L_0\subset \bigcup_{j=1}^{n} Y\setminus F_{i_j}$. But since $F$ is a chain and we chose finitely many indexes, there must exists some smallest element $F_{k} \in \{F_{i_1}, ..., F_{i_n}\}$. Note that $F_k\cap L_0$ is not empty and $F_k$ is disjoint from $\bigcup_{j=1}^{n} Y\setminus F_{i_j}$, a contradiction. Therefore $F$ is nonempty and belongs to $\mathcal{T}$. One verifies trivially that $F\preceq F_i$ for all $i\in I$ and so we can apply Zorn's Lemma Zorn's Lemma to obtain that $\mathcal{T}$ has an element $\overline F$ that is minimal for inclusion. Let $y_0$ be a point in $\overline F \cap L_0$. Note that, if $\omega(y_0)$ is the $\omega$-limit of $y_0$ by $T$, then it is also a compact invariant subset of $Y$, and since $u_0$ is recurrent, $\omega(y_0)$ must also intersect $L_0$. One deduces that $\omega(y_0)$ is also an element of $\mathcal{T}$ which is contained in $\overline F$, and thus is equal to $\overline F$ by minimality. But this implies $y_0\in\omega(y_0)$ and the claim is proved. 

%for all $i\in\N$, then taking $F=\bigcap_{i\in\N}F_i$, one observes trivially that $F$ is compact and invariant since $Y$ is compact. Also, note that $F\cap L_0=\bigcap_{i\in\N}\left(F_i\cap L_0\right)$ is nonempty as it is again a nested intersection of nonempty compact sets, since $L_0$ is also compact. This implies that $F\in\mathcal{T}$ and one checks that $F_i\preceq F$ for all positive $i$. We therefore apply Zorn's Lemma to obtain that $\mathcal{T}$ has an element $\overline F$ that is minimal for inclusion. Let $y_0$ be a point in $\overline F \cap L_0$. Note that, if $\omega(y_0)$ is the $\omega$-limit of $y_0$ by $T$, then it is also a compact invariant subset of $Y$, and since $u_0$ is recurrent, $\omega(y_0)$ must also intersect $L_0$. One deduces that $\omega(y_0)$ is also an element of $\mathcal{T}$ which is contained in $\overline F$, and thus is equal to $\overline F$ by minimality. But this implies $y_0\in\omega(y_0)$ and the claim is proved. 

Let then $Y^{\prime}$ be the closure of the forward orbit of $y_0$ by $T$,
which is a compact $T$-invariant set such that the restriction of $T$ to $Y^{\prime}$
 is both transitive and an extension of the shift $\sigma$. Let $%
\Lambda_0=\pi_1(Y^{\prime})$, a compact $g$-invariant set for which the
restriction of $g$ is also transitive. Given $p/q$ as in the statement,
consider 
$$
A_{p/q}=\left\{u\in\sigma_2;\, \hbox{for all } j\in {\rm Z}\!\!{\rm Z},
\left(\sum_{i= jq}^{j(q+1)-1} u_i\right)=p\right\}%
$$
which is invariant by $\sigma^{q}$, and let $Y_{p/q}=\pi_2^{-1}(A_{p/q})\cap
Y^{\prime}$, and $K_{p/q}=\pi_1(Y_{p/q})$. Note that the restriction of $%
\sigma^q$ to $A_{p/q}$ has strictly positive topological entropy (as it is
conjugated to the full shift on $\binom{{q}}{{p}}$ symbols). This implies
that the restriction of $T^q$ to $Y_{p/q}$ has positive topological entropy
and, since the cardinality of the fibers of $\pi_1$ is uniformly bounded,
the same holds for the restriction of $g^q$ to $K_{p/q}$. The second
assertion from the lemma follows directly from noticing that, if $y\in Y^{\prime}, x
=\pi_1(y)$ and $\widetilde x\in p^{-1}(x)\cap \widetilde \Lambda$, then $%
\widetilde g^{k}(\widetilde x)$ lies in $\widetilde
\Lambda+(\sum_{i=0}^{k-1}(\pi_2(y))_i, 0),$ concluding the proof. $\qed $

\vskip 0.2truecm

The following result is basically contained in subsection 6.1.2 of \cite
{contal}:

\begin{proposition}
\label{bigpropforcing} Let $A$ be a Birkhoff region of instability for some $%
f\in {\rm {Diff}_0^0(A)}$ with a lift $\widetilde{f}$ whose rotation set is
a non-degenerate interval. Then $f$ has a rotational topological horseshoe.
Moreover, for any nonempty open interval $J\subset \rho (\widetilde{f})$,
there exists a rotational topological horseshoe whose rotation set
intersects $J$.
\end{proposition}

Let us just explain how to perform the necessary modifications to that
subsection, so as to obtain this result. Using the same language and
definitions of \cite{contal}, the main idea in order to show that $f$ has a
rotational topological horseshoe is to apply Theorem M of \cite{lectal2}.
For this, consider an open interval $J\subset \rho( 
\widetilde{f})$ and choose some rational $s/r$ in $J$ such that $s/r$ is not
the rotation number for $\widetilde f$ of any point in the boundary of $A $.
Let $g=f^r$. Then there exists a {\it {maximal isotopy}} $I^{\prime}$
joining $g$ to the identity, that lifts to a maximal isotopy $\widetilde{I'}$
joining the identity in $\widetilde A$ to a lift $\widetilde g= \widetilde
f^{r}-(s,0)$ of $g$. This
can be done in such a way that the rotation set of $\widetilde g$ is an interval containing the origin, and such that the
rotation number of points in the upper boundary of $\widetilde A$ is not
null. For such maps, one can find a Brouwer-Le Calvez foliation $\mathcal{F}$
for $I^{\prime}$ that is lifted to a Brouwer-Le Calvez foliation $\widetilde{
\mathcal{F}}$ for $\widetilde{I'}$ and one can consider the set of
admissible $\widetilde{\mathcal{F}}$-transverse paths, as defined in \cite
{lectal1}. To show the existence of a rotational topological horseshoe, it
is then sufficient to show that there exists an $n$-admissible $\widetilde{ 
\mathcal{F}}$-transverse path $\widehat{\gamma'}$ such that $\widehat{\gamma'%
}$ has a $\widetilde{\mathcal{F}}$-transverse intersection with $\widehat{%
\gamma'}+(p,0)$ for some non-null integer $p$. This implies that there exists $n>0$ such that  $g$ has a
rotational topological horseshoe whose $\widetilde{g}$-rotation set contains 
$[0, p/n]$ if $p>0$, or contains $[p/n, 0]$ if $p$ is negative. In any case,
this implies that $f$ has a rotational topological horseshoe whose rotation
set contains an interval having $s/r$ as an endpoint, and thus intersecting $%
J$.

Most of the work contained in subsection 6.1.2 of \cite{contal} is concerned
with estimating the size of the rotation set contained in the rotational
topological horseshoe for $f$, and for such a reason a stronger hypothesis
than just asking $A$ to be a Birkhoff region of instability was assumed. But
in what concerns us here, which is to show that a rotational topological
horseshoe exists without requiring its rotation set to be of any specific
length, that extra hypothesis is unnecessary. The curve $\widehat{\gamma'}$
we need is obtained much in the same way as in the quoted subsection: One
shows that, assuming without loss of generality that the rotation number of
the upper boundary is strictly positive for $\widetilde g$, any point $%
\widetilde z$ in ${\rm I}\!{\rm R}{}\times\{1\}$ has a full $\widetilde{ 
\mathcal{F}}$-transverse trajectory that is equivalent to a $\widetilde{ 
\mathcal{F}}$-transverse simple curve $\widehat{\gamma}:{\rm I}\!{\rm R}%
{}\to \widetilde{A}$ satisfying $\widehat{\gamma}(t+1)=\widehat{\gamma}%
(t)+(1,0)$. In particular, using that $A$ is a Birkhoff region of
instability, there exist positive integers $N_0, N_1,$ a point $\widetilde
Z_{0,1}$ that is sufficiently close to the lower boundary of $\widetilde A$
with $\widetilde g^{N_0}(\widetilde Z_{0,1})$ sufficiently close to the
upper boundary of $\widetilde A$, so that its transverse path up to time $%
N_0 $ contains a subpath equivalent to $\widehat \gamma\mid_{[0,2]}$ which
starts at a leaf not intersected by $\widehat \gamma$, and another point $%
\widetilde Z_{1,0}$ that is sufficiently close to the upper boundary of $%
\widetilde A$ with $\widetilde g^{N_1}(\widetilde Z_{1, 0})$ sufficiently
close to the lower boundary of $\widetilde A$, and such that its transverse
path up to time $N_1 $ contains a subpath equivalent to $\widehat
\gamma\mid_{[0,2]},$ which ends at a leaf not intersected by $\widehat
\gamma $. The construction of $\widehat{\gamma'}$ then follows exactly as in
subsection 6.1.2 of \cite{contal}. $\qed$

A direct consequence of the two previous results is the following:

\begin{corollary}
\label{forcing} Let $A$ be a Birkhoff region of instability for some $f\in 
{\rm {Diff}_0^0(A)}$ and assume $\rho (\widetilde{f})$ is a non-degenerate
interval for some fixed lift $\widetilde{f}.$ Then there exists an open and
dense subset $E$ of $\rho (\widetilde{f})$ such that for any $p/q\in E$,
there exists:

\begin{itemize}
\item  Integers $s,r,L$ depending on $p/q$ and a rotational topological
horseshoe for $\widetilde{f}$ whose rotation set contains an interval $%
[s/r,(s+1)/r]$ with $s/r<p/q<p_1/q_1<(s+1)/r,$ where $p_1=p.r.L+1,q_1=q.r.L$;

\item  A compact and transitive $f$-invariant set $\Lambda ^{\prime
}=\Lambda ^{\prime }(p/q)$;

\item  Two compact $f$-invariant sets $G^0=G^0(p/q)$ and $G^1=G^1(p/q)$,
both contained in $\Lambda ^{\prime }$, such that the restriction of $f$ to
each of them has strictly positive topological entropy;

\item  Two compact sets $\widetilde{G^0},\widetilde{G^1}$ in $\widetilde{A}$
projecting surjectively onto $G^0$ and $G^1$ respectively, and such that $
\widetilde{f}^{q_1}(\widetilde{G^0})-(p.r.L,0)=\widetilde{G^0}$ and $
\widetilde{f}^{q_1}(\widetilde{G^1})-(p_1,0)=\widetilde{G^1}$;
\end{itemize}
\end{corollary}

{\it Proof.} Let $E$ be the union of the interior of the rotation sets of
all rotational topological horseshoes of $f$. Proposition~\ref
{bigpropforcing} shows that $E$ is dense, and as it is the union of open
intervals, it is also open. If $p/q$ is a rational point in $E$, then by
definition one can find a rotational horseshoe whose rotation set contains $%
p/q$ in its interior, and let then $s, r$ be integers such that $%
s/r<p/q<(s+1)/r$ and such that $[s/r,(s+1)/r]$ is also contained in the
interior of the rotation set of this topological horseshoe. If we then
choose $L$ sufficiently large, than taking $p_1=p.r.L+1,q_1=q.r.L$, we get
that $p_1/q_1<(s+1)/r$. This shows the first item.

% such that any rational in its interior is also in the interior of the rotation set of a rotational topological horseshoe follows directly from Proposition~\ref{bigpropforcing}, as we see that for any open set in $\rho (\widetilde{f})$ there exists a rotational topological horseshoe whose rotation set intersects this set. By the choice of $E$, given $p/q\in E$, one can of course find a horseshoe satisfying the first item, that is, the $\widetilde{f}$-rotation set of the horseshoe contains some interval $[s/r,(s+1)/r],$ and $p/q,p_1/q_1$ both belong to its interior (for some $L$ as in the statement of the corollary).

To get the other itens of the proposition, note first that, by setting $%
p^*=rp-sq, p_1^*=\frac{1}{r}\cdot(rp_1-sq_1)=p_1-sqL$, then $p^*/q=
r(p/q)-s,\, p_1^*/(Lq)=r(p_1/q_1)-s$ and so $0<p^*/q<p_1^*/(Lq)<1$. We apply
Lemma ~\ref{transitiveentropy} with $g=f^r,\widetilde{g}=\widetilde{f}%
^r-(s,0)$, giving us a set $\Lambda_0$ in $A$ which is transitive for $g$,
as well as two sets $K_{p^*/q}, K_{p_1^*/(Lq)}$, both contained in $%
\Lambda_0 $, such that the former is invariant by $g^q$ and the later is
invariant by $g^{Lq}$, and such that the topological entropy of the
restriction of $g^{Lq}$ to either $K_{p^*/q}$ or $K_{p_1^*/(Lq)}$ is strictly
positive. Furthermore, there exists sets $\widetilde{K}_{p^*/q}, \widetilde{K%
}_{p_1^*/(Lq)}$ in $\widetilde{A}$ projecting onto $K_{p^*/q},
K_{p_1^*/(Lq)} $ respectively and such that $\widetilde g^q(\widetilde{K}%
_{p^*/q})= \widetilde{K}_{p^*/q}+(p^*,0)$ and such that $\widetilde
g^{(Lq)}( \widetilde{K}_{p_1^*/(Lq)})= \widetilde{K}_{p_1^*/(Lq)}+(p_1^*,0)$.

Define now $\Lambda ^{\prime }=\bigcup_{i=0}^{r-1}f^i(\Lambda _0)$, and note
that, as $\Lambda_0$ was an invariant set for $g=f^r$ and the restriction of 
$g$ to this set was transitive, then $\Lambda^{\prime}$ is invariant for $f$
and the restriction of $f$ to it is also transitive. This gives us the
second item of the corollary. Define also 
$$
G^0=\bigcup _{i=0}^{r-1}f^i\left( \bigcup
_{j=0}^{q-1}g^j(K_{p^*/q})\right),\, \hbox{and }G^1=\bigcup
_{i=0}^{r-1}f^i\left( \bigcup _{j=0}^{q_1-1}g^j(K_{p_1^*/(Lq)})\right),%
$$
which are both subsets of $\Lambda^{\prime}$ as they are each contained in
the $f$-orbit of $K_{p^*/q}$ and $K_{p_1^*/(Lq)}$, respectively, and $%
K_{p^*/q}$ and $K_{p_1^*/(Lq)}$ are contained in $\Lambda^{\prime}$ which is
invariant. One also observes that $G^0$ and $G^1$ are $f$ invariant by
construction, and the restriction of $f$ to both this sets has strictly
positive topological entropy. This gives us the third item of the
proposition.

Finally, set 
$$
\widetilde{G^0}=\bigcup_{i=0}^{r-1}\widetilde{f}^i\left( \bigcup
_{j=0}^{q-1} \widetilde{g}^j( \widetilde{K_{p^*/q}})\right),\, \hbox{and } 
\widetilde{G^1}=\bigcup_{i=0}^{r-1}\widetilde{f}^i\left( \bigcup
_{j=0}^{q_1-1}\widetilde{g}^j( \widetilde{K_{p_1^*/(Lq)}})\right).%
$$
Note that, as $\widetilde f$ and $\widetilde g$ commute, we get that

\begin{align*}
\widetilde{f}^{rq}(\widetilde{G^0})&=\widetilde g^q(\widetilde{G^0})+(sq,0)\\
&=\bigcup_{i=0}^{r-1}\widetilde{f}^i\left( \bigcup _{j=0}^{q-1}\widetilde{g}^j( 
\widetilde g^q(\widetilde{K_{p^*/q}}))\right)+(sq,0)\\
&=\bigcup_{i=0}^{r-1}\widetilde{f}^i\left( \bigcup _{j=0}^{q-1}\widetilde{g}^j( 
\widetilde{K_{p^*/q}})\right)+(sq,0)+(p^*,0)\\
&= \widetilde{G^0}+(rp,0)
\end{align*}

and so $\widetilde{f}^{q_1}(\widetilde{G^0})=\widetilde{f}^{L\cdot rq}( 
\widetilde{G^0})=\widetilde{G^0}+(L\cdot rp,0)$. A similar computation shows
that $\widetilde{f}^{q_1}(\widetilde{G^1})= \widetilde{G^1}+(sqL+p^*_1,0)= 
\widetilde{G^1}+(p_1,0) $ ending the claim $\qed$ 

We remark that from the
corolary, every point in the set $G^0(p/q)$ has rotation number $p/q$ for $%
\widetilde f$, and every point in $G^1(p/q)$ has rotation number $%
p_1/q_1$ for $\widetilde f$.

\subsection{On maximal invariant sets}

Let $\gamma \subset S^1\times ]0,1[$ be a homotopically non-trivial simple
closed curve and, as we defined before, let $\gamma ^{-}$ be the connected
component of $\gamma ^c$ which contains $S^1\times \{0\}$ (similary for $%
\gamma ^{+}$ and $S^1\times \{1\}$).

If we consider the sets

$$
B_{0,\gamma }^s={\bigcap }_{n\leq 0} f^n(\overline{\gamma ^{-}}), \text{ } B_{0,\gamma
}^u={\bigcap }_{n\geq 0} f^n(\overline{\gamma ^{-}}), 
$$
$$
B_{1,\gamma }^s={\bigcap }_{n\leq 0} f^n(\overline{\gamma ^{+}}) \text{ and }
B_{1,\gamma }^u={\bigcap }_{n\geq 0} f^n(\overline{\gamma ^{+}}), 
$$
we get that,

$$
f(B_{0,\gamma }^s)\subset B_{0,\gamma }^s,\text{ }f^{-1}(B_{0,\gamma
}^u)\subset B_{0,\gamma }^u, 
$$

$$
f(B_{1,\gamma }^s)\subset B_{1,\gamma }^s\text{ and }f^{-1}(B_{0,\gamma
}^u)\subset B_{0,\gamma }^u. 
$$

Denote by $\widehat{B}_{0,\gamma }^s$ the connected component of 
$B_{0,\gamma }^s$ which contains $S^1\times \{0\}$ and define similarly $\widehat{B}_{0,\gamma }^u$. 
Analogously, let $\widehat{B}_{1,\gamma }^s$ be the connected component of $B_{1,\gamma }^s$
which contains $S^1\times \{1\}$ and define similarly $\widehat{B}_{1,\gamma }^u$. 

The next result appears in Le Calvez \cite{inst} and even in Birkhoff's
paper \cite{birk}.

\begin{lemma}
\label{bmenos} Let $f:A\rightarrow A$ be an orientation and boundary
components preserving homeomorphism, which has the curve intersection
property. Then, for any $\gamma $ as above, $\widehat{B}_{0,\gamma }^s, 
\widehat{B}_{1,\gamma }^s,\widehat{B}_{0,\gamma }^u$ and $\widehat{B}%
_{1,\gamma }^u$ intersect $\gamma $$.$
\end{lemma}

As $f(\widehat{B}_{0,\gamma }^s)\subset \widehat{B}_{0,\gamma }^s,$ $f( 
\widehat{B}_{1,\gamma }^s)\subset \widehat{B}_{1,\gamma }^s,$ $f^{-1}( 
\widehat{B}_{0,\gamma }^u)\subset \widehat{B}_{0,\gamma }^u,$ and $f^{-1}( 
\widehat{B}_{1,\gamma }^u)\subset \widehat{B}_{1,\gamma }^u,$ we get that
the maximal invariant sets $\xi ^{-}(\gamma )$ and $\xi ^{+}(\gamma )$
satisfy the following conditions:

$$
{\cap }_{n\geq 0}f^{-n}(\widehat{B}_{0,\gamma }^u)={\cap }_{n\geq 0}f^n( 
\widehat{B}_{0,\gamma }^s)=\xi ^{-}(\gamma ) 
$$

$$
\text{ and } 
$$

$$
{\cap }_{n\geq 0}f^{-n}(\widehat{B}_{1,\gamma }^u)={\cap }_{n\geq 0}f^n( 
\widehat{B}_{1,\gamma }^s)=\xi ^{+}(\gamma ) 
$$

\section{Proofs}

\subsection{Proof of Theorem 2}

Let us remember our set of hypotheses:\ $f\in {\rm {Diff}_0^{1+\varepsilon
}(A)}$ for some $\varepsilon >0,$ $A$ is a Birkhoff region of instability
and for some fixed lift $\widetilde{f},$ $\rho (\widetilde{f})$ has interior.

Let $E$ be the set from Corollary~\ref{forcing}, which we further assume
does not contain the rotation number of the boundaries. Fix some $p/q\in E$
and consider, as in Corollary~\ref{forcing}, the sets $\Lambda ^{\prime
},G^0(p/q),G^1(p/q)$ as well as the integers $r,s,L,p_1$ and $q_1$, all
depending on $p/q$. As $h(${\it f}$\mid _{G^0})>0$ and $h(${\it f}$\mid _{G^1})>0,$ there
exist two hyperbolic ergodic Borel non-atomic $f$-invariant measures $\mu _{p/q}$
and $\mu _{p_1/q_1},$ such that 
${\rm {supp}}(\mu_{p/q})\subset G^0(p/q)$ and 
${\rm {supp}}(\mu _{p_1/q_1})\subset G^1(p/q)$. Note that, by the choice of 
$G^0(p/q)$ (respectively $G^1(p/q)$), every point  in it has rotation number $p/q$ (respectively $p_1/q_1$). 
From Lemma \ref{andre2}, pick inaccessible points $z_0\in {\rm {%
supp}(}\mu _{p/q}{\rm )}$ and $z_1\in {\rm {supp}(}\mu _{p_1/q_1}{\rm ).}$

Lemma \ref{andre1} implies that there are four hyperbolic periodic saddle points, 
$y_0$ and $y_0^{\prime}$, $y_1$ and $y_1^{\prime },$ such that $z_0$ is enclosed 
by the rectangle
determined by compact subarcs of stable and unstable branches at $y_0$ and
at $y_0^{\prime },$ an analogous statement holding for $z_1$ and $y_1$ and $%
y_1^{\prime }.$ 
%(the above mentioned branches at $Q_0$ are denoted $\lambda^s_0$ and $\lambda^u_0$).
The corners of each rectangle are either the saddles, or $C^1$-transverse
intersections between stable branches at one saddle and unstable branches at
the other. In particular, the $\lambda $-lemma implies that for each of
these four periodic points $y_0,y_0^{\prime },y_1$ and $y_1^{\prime },$
there are $C^1$-transverse homoclinic intersections.

As the rectangles can be chosen in an arbitrarily small way, the rotation
numbers of $y_0$ and $y_0^{\prime }$ are both equal to $p/q$ (but their
periods might be larger than $q$) and the rotation numbers of $y_1$ and 
$y_1^{\prime }$ are both equal to $p_1/q_1,$ which is different from $p/q$.
So, as $\overline{W^u(y_0)}=\overline{W^u(y_0^{\prime })}$, $\overline{%
W^u(y_1)}=\overline{W^u(y_1^{\prime })}$, $\overline{W^s(y_0)}=\overline{%
W^s(y_0^{\prime })}$ and $\overline{W^s(y_1)}=\overline{W^s(y_1^{\prime })}$%
, %$ and $\overline{W^s(Q_{0\text{ or }1})}=\overline{W^s(Q_{0 \text{ or }%
%1}^{\prime })}$
%
%
%$
%\overline{W^u(Q_{0\text{ or }1})}=\overline{W^u(Q_{0\text{ or }1}^{\prime })}
%$ and $\overline{W^s(Q_{0\text{ or }1})}=\overline{W^s(Q_{0 \text{ or }%
%1}^{\prime })}$
and there is a dense orbit in $\Lambda ^{\prime }$, arbitrarily large
positive iterates of the interior of the rectangle enclosing $z_0$ intersect
the interior of the rectangle enclosing $z_1$ and vice-versa. And this
implies (see Lemma 24 of \cite{c1epsilon}) that for some integer $i,$ there
exists an unstable branch at $y_0$ that has a topologically transverse
intersection with a stable branch at $f^i(y_1)$ and an unstable branch at $%
f^i(y_1)$ has a topologically transverse intersection with a stable branch
at $y_0.$ As the rotation number of $y_0$ is $p/q$ and the rotation number
of $f^i(y_1)$ is not $p/q,$ the theorem follows from the $C^0$-$\lambda $%
-lemma that holds for topologically transverse intersections, 
see Proposition \ref{toptrans}. $\qed $

\subsection{Proof of Lemma~1}

Under the lemma hypotheses, Theorem~\ref{dois} implies that for any rational
point $p/q$ in the set $E,$ we can find an hyperbolic periodic saddle $%
z_{p/q}$ with rotation number $p/q$ and unstable and stable branches, $%
\lambda _{p/q}^u$ and $\lambda _{p/q}^s,$ both at $z_{p/q},$ which intersect
at a point $w_{p/q},$ such that if we concatenate the arc in $\lambda
_{p/q}^u$ from $z_{p/q}$ to $w_{p/q}$ to the arc in $\lambda _{p/q}^s$ from $%
w_{p/q}$ to $z_{p/q}$, then we get a homotopically non-trivial closed curve $%
\gamma _{p/q}$ contained in the interior of $A$ (because the rotation
numbers on the boundary components do not lie in $E$).

\begin{proposition}
\label{acumul} For each $p/q\in E,$ $\overline{\lambda _{p/q}^u}$ and $
\overline{\lambda _{p/q}^s}$ intersect $S^1\times \{0\}$ and $S^1\times
\{1\}.$
\end{proposition}

{\it Proof. }Fixed some $p/q \in E,$ let $q.k_{p/q}$ be twice the period of
$z_{p/q}$. As $A$ is a Birkhoff region of instability, $\cup _{n\geq
0}f^{n.q.k_{p/q}}(\gamma_{p/q})$ accumulates on both boundary components, $%
S^1\times \{0\}$ and $S^1\times \{1\}.$ An analogous statement holds for $%
\cup _{n\geq 0}f^{-n.q.k_{p/q}}(\gamma_{p/q}).$ The proposition follows from
the $f^{q.k_{p/q}}$-positive invariance of the arc in $\lambda_{p/q}^s$ from 
$z_{p/q}$ to $w_{p/q}$ and the $f^{q.k_{p/q}}$-negative invariance of the
arc in $\lambda_{p/q}^u$ from $z_{p/q}$ to $w_{p/q}.$ $\qed $

\vskip 0.2truecm

The curve $\gamma _{p/q}$ may not be simple, but in any case, $(\gamma
_{p/q})^c$ still has one connected component which contains $S^1\times
\{0\}, $ denoted $\gamma _{p/q}^{-},$ another connected component which
contains $S^1\times \{1\},$ denoted $\gamma _{p/q}^{+}$ and maybe other
contractible components.

\begin{proposition}
For any $p/q\in E$ and any integer $n>1$, 
$$
\xi _{1/n}^{-}\subset \overline{\gamma _{p/q}^{-}}\text{ and }\xi
_{1/n}^{+}\subset \overline{\gamma _{p/q}^{+}}. 
$$
In particular, there exist $f$-invariant continua $K^{-}\supset S^1\times
\{0\}$ and $K^{+}\supset S^1\times \{1\}$ such that $\xi _{1/n}^{-}=K^{-}$
and $\xi $$_{1/n}^{+}=K^{+}$ for all sufficiently large $n$.
\end{proposition}

{\it Proof.} First, note that Proposition \ref{acumul} implies that $\xi
_{1/n}^{-}$ can not contain $\lambda _{p/q}^u$ or $\lambda _{p/q}^s$ because
both branches accumulate on $S^1\times \{1\}.$ Moreover, if for some $n>1,$ $%
\xi $$_{1/n}^{-}$ is not contained in $\overline{\gamma _{p/q}^{-}},$ then
as $\xi _{1/n}^{-}$ is $f$-invariant, connected and contains $S^1\times
\{0\} $, $\xi _{1/n}^{-}$ would contain two sequences of points, one in $%
\gamma _{p/q}^{-}$ and one in $\left( \overline{\gamma _{p/q}^{-}}\right) ^c,
$ both converging to $z_{p/q}.$ At least one of them converges to $z_{p/q}$
through the local quadrant at $z_{p/q}$ adjacent to $\lambda _{p/q}^u$ and $%
\lambda _{p/q}^s.$ So, Proposition 6, item 2 of \cite{bounded} (which says
that any $f$-invariant continuum which contains an hyperbolic saddle
periodic point $p$ and accumulates on $p$ through a certain local quadrant $Q
$ at $p,$ must contain at least one branch at $p$ adjacent to $Q)$ implies
that $\xi _{1/n}^{-}$ contains either $\overline{\lambda _{p/q}^u}$ or $
\overline{\lambda _{p/q}^s},$ a contradiction as explained above. 
%because they both intersect 
%$S^1\times \{1\}.$ 
One shows by a similar argument that $\xi $$_{1/n}^{+}$ is contained in $
\overline{\gamma _{p/q}^{+}}$.

The above argument implies that, for any homotopically non-trivial simple
closed curve $\gamma $ contained in the interior of $A,$ such that $\gamma
_{p/q}^{-}\subset \gamma ^{-}$, $\xi $$^{-}(\gamma )$ must be equal to the
connected component of the maximal invariant set contained in the closure of 
$\gamma _{p/q}^{-}$ which contains $S^1\times \{0\}$. A similar statement
holds for $\xi $$^{+}(\gamma )$. Therefore, if $N>0$ is such that $S^1\times
\{1-1/N\}\subset \gamma _{p/q}^{+}$ and $S^1\times \{1/N\}\subset \gamma
_{p/q}^{-},$ then for all $n\geq N,$ $\xi $$_{1/n}^{-}=\xi _{1/N}^{-}:=K^{-}$
and $\xi $$_{1/n}^{+}=\xi _{1/N}^{+}:=K^{+}$. $\qed$

\vskip 0.2truecm

Note that the intersection between $K^{-}$ and $K^{+}$ must be empty because
otherwise one could apply as before, Proposition 6, item 2 of \cite{bounded}
and get that either $K^{-}$ or $K^{+}$ (maybe both), contains $\overline{%
\lambda _{p/q}^u}$ or $\overline{\lambda _{p/q}^s}.$ As explained above,
this is a contradiction.

So, as $(K^{-})^c$ and $(K^{+})^c$ are both connected, $(K^{-}\cup
K^{+})^c:=A^{*}$ is an open $f$-invariant essential annulus contained in $A.$
As $A$ is a Birkhoff region of instability, $\partial A^{*}$ intersects both 
$S^1\times \{0\}$ and $S^1\times \{1\}.$

We finish with the following proposition.

\begin{proposition}
\label{propaux} For all rationals $p/q\in E$, with the exception of at most
two points, $\gamma _{p/q}$ is contained in $A^{*}$.
\end{proposition}

{\it Proof. }If, for some $p/q$ in $E$, $\gamma _{p/q}$ intersects $K^{-},$
then as $K^{-}$ is compact and $f$-invariant, $z_{p/q}\in K^{-},$ and it is
clearly accessible from $(K^{-})^c.$ This happens because $K^{-}\subset 
\overline{\gamma _{p/q}^{-}}$ and every point in the upper boundary of $
\overline{\gamma _{p/q}^{-}}$ (which contains $z_{p/q}$) is the endpoint of
a $C^1$ endcut contained in the complement of $\overline{\gamma _{p/q}^{-}}.$
This is a trivial consequence of the fact that $\gamma _{p/q}$ is a
piecewise $C^1$ curve. So, Lemma \ref{acessperiod} implies that all
accessible periodic points in $K^{-}$ have $p/q$ as rotation number (the
rotation number of the prime ends compactification of $(K^{-})^c$ is equal
to $p/q$). A similar argument holds for $K^{+},$ and so, with the exception
of at most two rational numbers in $E$, for all other rationals $p/q\in E$, $%
\gamma _{p/q}$ avoids $K^{-}\cup K^{+},$ in other words, $\gamma _{p/q}$ is
contained in $A^{*}.$ $\qed $

\vskip0.2truecm

Since $E$ is an infinite set, this concludes the proof of the lemma. $\qed$

\subsection{Proof of Theorem 1}

From Proposition~\ref{propaux}, there exists $p/q\in E$ such that $\gamma
_{p/q}\subset A^{*}.$ Moreover, $\gamma _{p/q}=z_{{p/q}}\cup \lambda
_{\rm{comp}}^u\cup \lambda _{\rm{comp}}^s,$ where $\lambda _{\rm{comp}}^u$ is the closed
arc in $\lambda _{p/q}^u$ from $z_{p/q}$ to $w_{p/q}$ and $\lambda _{\rm{comp}}^s$
is the closed arc in $\lambda _{p/q}^s$ from $w_{p/q}$ to $z_{p/q}$, where $%
w_{p/q}$ is a point in the intersection between $\lambda _{p/q}^u$ and $%
\lambda _{p/q}^s$.

Denote the lower (resp. upper) connected component of the boundary of $A^{*}$
by $K^{*0}\subset K^{-}$ (resp. $K^{*1}\subset K^{+}$). For $\gamma \subset
A $ a homotopically non-trivial simple closed curve, which is also contained
in $A^{*},$ we can define $\xi $$^{*-}(\gamma )$ as the connected component
of the maximal invariant set contained in $\overline{\gamma ^{-}}\cap \overline{A^{*}}$
which contains $K^{*0}.$ Analogously for $\xi $$^{*+}(\gamma ).$ From the
construction of $K^{-}$ and $K^{+}$ in Lemma \ref{umaux}, we get that 
\begin{equation}
\label{umapouca}\xi ^{*-}(\gamma )=K^{*0}\text{ and }\xi ^{*+}(\gamma
)=K^{*1}. 
\end{equation}
Now consider homotopically non-trivial simple closed curves $\alpha _1$ and $%
\alpha _0,$ both contained in $A^{*},$ such that $\alpha _1\subset \gamma
_{p/q}^{+}$ and $\alpha _0\subset \gamma _{p/q}^{-}.$

Also consider the sets $\widehat{B}_{0,\alpha _1}^{*s},$ $\widehat{B}%
_{0,\alpha _1}^{*u},$ $\widehat{B}_{1,\alpha _0}^{*s}$ and $\widehat{B}%
_{1,\alpha _0}^{*u}$ defined in Subsection 2.5, but now with respect to 
 $\overline{A^{*}}.$

Clearly,

\begin{equation}
\label{duaspouca} 
\begin{array}{c}
{\cap }_{n\geq 0} f^n(\widehat{B}_{0,\alpha _1}^{*s})= {\cap }_{n\geq 0}
f^{-n}(\widehat{B}_{0,\alpha _1}^{*u})=K^{*0} \\ \text{ and } \\ {\cap }%
_{n\geq 0} f^n(\widehat{B}_{1,\alpha _0}^{*s})= {\cap }_{n\geq 0} f^{-n}( 
\widehat{B}_{1,\alpha _0}^{*u})=K^{*1}. 
\end{array}
\end{equation}

%Proposition \ref{acumul} implies that there exist compact arcs, $\alpha_{p/q}^u$ contained in $\lambda _{p/q}^u$ and $\alpha _{p/q}^s$ contained in $\lambda _{p/q}^s,$ such that both $\alpha _{p/q}^u$ and $\alpha _{p/q}^s$intersect $\alpha _0$ and $\alpha _1$.

The above equalities imply that for any $z\in \widehat{B}_{0,\alpha _1}^{*s}$%
, its $\omega $-limit set is contained in $K^{*0}$ and for any $w\in 
\widehat{B}_{0,\alpha _1}^{*u},$ its $\alpha $-limit set is also contained
in $K^{*0}.$ And analogously, %if $z \in \widehat{B}_{0,\alpha _1}^{*u}$, 
%its backward orbit converges to $K^*_0$, 
for any $z\in \widehat{B}_{1,\alpha _0}^{*s}$, its $\omega $-limit set is
contained in $K^{*1}$ and for any $w\in \widehat{B}_{0,\alpha _1}^{*u}$, its 
$\alpha $-limit set is also contained in $K^{*1}.$ %
%and if $z \in \widehat{B}_{1,\alpha _0}^{*u}$, 
%its backward orbit converges to $K^*_1$. 
Since $K^{*0}$ and $K^{*1}$ are, respectively, the lower and the upper
connected components of the boundary of $A^{*}$, we will conclude the proof
of the theorem by showing that for some integer $n\geq 0$, $f^n(\widehat{B}%
_{0,\alpha _1}^{*u})\cap f^{-n}(\widehat{B}_{1,\alpha _0}^{*s})$ and $f^n(
\widehat{B}_{1,\alpha _0}^{*u})\cap f^{-n}(\widehat{B}_{0,\alpha _1}^{*s})$
are both non-empty.

\begin{lemma}
The sets $\widehat{B}_{0,\alpha _1}^{*s}$ and $\widehat{B}_{1,\alpha _0}^{*s}
$ intersect $\lambda _{{p/q}}^u$ in a topologically transverse way, and
analogously, $\widehat{B}_{0,\alpha _1}^{*u}$ and $\widehat{B}_{1,\alpha
_0}^{*u}$ intersect $\lambda _{{p/q}}^s$ also in a topologically transverse
way.
\end{lemma}

{\it Proof. }The proof is analogous in all four cases, for $\widehat{B}%
_{0,\alpha _1}^{*s},\widehat{B}_{1,\alpha _0}^{*s},\widehat{B}_{0,\alpha
_1}^{*u}$ and $\widehat{B}_{1,\alpha _0}^{*u}.$ So, without loss of
generality, let us only analyze $\widehat{B}_{0,\alpha _1}^{*s}.$

Let $\Theta $ be a connected component of the intersection between
%$$
$\widehat{B}_{0,\alpha _1}^{*s}$ %\cap 
and the closed annulus bounded by $\alpha _0$ and $\alpha _1$ 
%$$
which intersects $\alpha _1$ (see Lemma {\ref{bmenos}). As $\Theta$ intersects 
$\alpha _1$, as $\widehat{B}_{0,\alpha _1}^{*s}$ is connected and contains $K^{*0}$, 
and as each connected component of $\widehat{B}_{0,\alpha _1}^{*s}$ is entirely 
contained in the closed annulus bounded by $\alpha_1$ and the lower boundary of $A$, 
one obtains that $\Theta$ must also 
intersect $\alpha _0$. 

From expression (2), $\widehat{B}%
_{0,\alpha _1}^{*s}$ does not intersect $\lambda _{p/q}^s.$ Assume, for a
contradiction, that $\Theta $ does not intersect $\lambda _{{p/q}}^u$ in a
topologically transverse way. Let $\varepsilon $$>0$ be sufficiently small in a
way that the $\varepsilon $-neighborhood of $\lambda _{\rm{comp}}^u,$ denoted $V,$
is contractible in $A^{*}$ and disjoint from $\alpha _0$ and $\alpha _1.$
Then, as the endpoints of $\lambda _{\rm{comp}}^u$ do not belong to $\widehat{B}%
_{0,\alpha _1}^{*s},$ by Lemma \ref{l.pertseparada}, one can find a curve $%
\mu _\varepsilon $ contained in $V$, with the same endpoints as $\lambda
_{\rm{comp}}^u,$ which is disjoint from $\Theta $. But as $V$ is contractible, $%
\mu _\varepsilon\cup \lambda _{\rm{comp}}^s$ is an homotopically non-trivial closed curve
separating $\alpha _0$ and $\alpha _1$ and disjoint from $\Theta .$ And this
is a contradiction since $\Theta $ is connected and intersects $\alpha _0$
and $\alpha _1$. \qed
}

\vskip0.2truecm

Thus, from the above lemma and Proposition \ref{toptrans}, $f^n(\widehat{B}%
_{0,\alpha _1}^{*u})$ and $f^n(\widehat{B}_{1,\alpha _0}^{*u})$ both contain
subcontinua which accumulate on compact sub arcs of $W^u(z_{p/q})$ in the
Hausdorff topology as $n\rightarrow \infty ,$ and analogously, $f^{-n}( 
\widehat{B}_{0,\alpha _1}^{*s})$ and $f^{-n}(\widehat{B}_{1,\alpha _0}^{*s})$
both contain subcontinua which accumulate on compact sub arcs of $%
W^s(z_{p/q})$ in the Hausdorff topology as $n\rightarrow \infty $. And this
means that, %As $\lambda _{i_0}^u$ and $\lambda _{i_0}^s$ have
%topologically transverse intersections, 
for a sufficiently large $n>0,$ 
$$
\begin{array}{c}
f^n( 
\widehat{B}_{0,\alpha _1}^{*u})\text{ intersects both }f^{-n}(\widehat{B}%
_{0,\alpha _1}^{*s})\text{ and }f^{-n}(\widehat{B}_{1,\alpha _0}^{*s}), \\ 
\text{and} \\ f^n(\widehat{B}_{1,\alpha _0}^{*u})\text{ intersects both }%
f^{-n}(\widehat{B}_{0,\alpha _1}^{*s})\text{ and }f^{-n}(\widehat{B}%
_{1,\alpha _0}^{*s}). 
\end{array}
$$
Denoting a point in $f^n(\widehat{B}_{0,\alpha _1}^{*u})\cap f^{-n}(\widehat{%
B}_{1,\alpha _0}^{*s})$ by $z^{+}$ and a point in $f^n(\widehat{B}_{1,\alpha
_0}^{*u})\cap f^{-n}(\widehat{B}_{0,\alpha _1}^{*s})$ by $z^{-},$ the
theorem is proved. $\qed $

\subsection{Proof of Theorem 3}

If the rotation set of $\widetilde f$ is the singleton $\{a\}$, there is
nothing to be done, since it follows easily that each point in the boundary
of $A$ has rotation number $a$ (and one can even show that every point in $A$
has rotation number $a$). So we can assume that the rotation set of $%
\widetilde f$ has nonempty interior and therefore we are in the hypotheses
of Theorem~\ref{um}. Let $A^*$ be given by this result, which is obtained as
in Lemma~1. Note that, from Theorem~2 we know that there exists $E$ which is
open and dense in $\rho(\widetilde f)$ such that for any rational $p/q$ in $%
E $ we find the homotopically non-trivial closed curve $\gamma_{p/q}$ as
described there, and from the end of the proof of Lemma~1, we know that all
but at most two of these curves are contained in $A^*$.

As is done for the disk, one can consider the prime ends compactification of 
$A^*,$ by adding two circles in order to obtain $A^{\prime}= A^*\sqcup
S^1\sqcup S^1,$ so that $A^{\prime}$ is homeomorphic to $A$, and the
restriction of $f$ to $A^{*}$ extends continuously to an homeomorphism $h$
of $A^{\prime}$, with a lift $\widetilde{h}$ to the universal covering of $%
A^{\prime}.$ The rotation number of the lower boundary component of $%
A^{\prime}$ is the prime ends rotation number of $K^{*0}$ and the rotation
number of the upper boundary component of $A^{\prime}$ is the prime ends
rotation number of $K^{*1}.$ Note that $A^{\prime}$ is a Mather region of
instability for $h$. Note also that, since every rational in $E$ is the $
\widetilde{h}$-rotation number of a point in $\overline{A^{*}}$, then the
rotation set of $\widetilde{h}$ must be the same as that of $\widetilde f$,
since rotation sets are closed. Now, Theorem C of \cite{contal} shows that
every point in the rotation set of $\widetilde{h}$ is realized by a compact $%
h$-invariant set, which implies that, except maybe for the two prime ends
rotation numbers of the boundary components of $A^{*}$, every point in the
rotation set of $\widetilde{f}$ is realized by a compact $f$-invariant
subset in $A^{*}.$ $\qed$

\vskip0.2truecm

%{\it Acknowledgements: }

%\centerline{\bf Figure captions.}

%\begin{itemize}
%\item[Figure 1. ]  Diagram showing $\widehat{A}.$

%\item[Figure 2. ]  Diagram showing the set $\Gamma _N.$

%\item[Figure 3. ]  Diagram showing that $\left[ \Gamma _2\cap V_a^{-}\right]
%\subset \Gamma _{1a,down}\cup \Gamma _{1a,up}.$

%\item[Figure 4. ]  Diagram showing that either $\left[ \Gamma _2\cap
%V_a^{-}\right] \subset \Gamma _{1a,down}$ or $\left[ \Gamma _2\cap
%V_a^{-}\right] \subset \Gamma _{1a,up}.$

%\item[Figure 5. ]  Diagram showing the sets $\Gamma $, $\Gamma -(1,0)$ and $%
%\gamma .$

%\item[Figure 6. ]  Diagram showing the sets $\overline{p(\Gamma )}\subset
%\gamma _E^{-}.$
%\end{itemize}

%\begin{center}
%\mbox{\includegraphics[width=11cm]{fig1.ps}}
%\end{center}

\end{document}